

\baselineskip=14pt
\parskip=10pt

\font\eightrm=cmr8 
\font\eighttt=cmtt8
\magnification=\magstephalf

\def\1{{\overline{1}}}
\def\2{{\overline{2}}}
\parindent=0pt
\overfullrule=0in

\def\frac#1#2{{#1 \over #2}}
\bf
\centerline
{
Automated Counting of Towers (\`A La Bordelaise)
}
\centerline
{
[Or: Footnote to p. 81 of the Flajolet-Sedgewick Chef-d'oevre]
}

\rm
\bigskip
\centerline{ {\it Shalosh B.   EKHAD} and
{\it Doron  ZEILBERGER}\footnote{$^1$}
{\eightrm  \raggedright
Department of Mathematics, Rutgers University (New Brunswick),
Hill Center-Busch Campus, 110 Frelinghuysen Rd., Piscataway,
NJ 08854-8019, USA.
{\eighttt zeilberg  at math dot rutgers dot edu} ,
\hfill \break
{\eighttt http://www.math.rutgers.edu/\~{}zeilberg/} .
Written: 12.12.12.
Accompanied by the Maple package  TOWERS available from
{\eighttt http://www.math.rutgers.edu/\~{}zeilberg/mamarim/mamarimhtml/migdal.html} ,
where one can also find lots of output that should be considered part of this article.
Supported in part by the NSF.
}
}

\bigskip
{\it \qquad  \qquad  \qquad  \qquad  \qquad  \qquad  \qquad  \qquad  
\qquad  \qquad   \qquad  \qquad  
\`A la m\'emoire de Philippe FLAJOLET}

{\bf Introduction}

One of our favorite theorems in enumerative combinatorics, 
whose proof-from-the-book by Jean B\'etr\'ema and Jean-Guy Penaud[BeP] is
succinctly outlined on p. 81 of  the Flajolet-Sedgewick[FS] {\it bible},
(that, in turn, is based on Mireille Bousquet-M\'elou's ``insightful presentation'' [Bo]), is
the $3^n$ theorem[GV] of Dominique Gouyou-Beauchamps  and Xavier Viennot.

This amazing result  asserts that there are {\bf exactly} $3^n$ ways of forming a tower of $n+1$ domino pieces
such that the bottom floor consists of one or more {\it contiguous} pieces,
and every piece at a higher floor touches one or two pieces on the floor right under
it in such a way that the common boundary is exactly (the left- or right-) half of either piece.

For example, here are $100$ such towers with $35$ pieces:

{\tt http://www.math.rutgers.edu/\~{}zeilberg/viennot/XR1b.html} \quad ,

and here are {\bf all} $27$ such towers with four pieces:

{\tt http://www.math.rutgers.edu/\~{}zeilberg/viennot/X4.html} \quad .

We love this proof (and theorem) {\bf so much} that we wrote the admiring {\it expos\'e} [Z1].

{\bf Why the present article?}

Physcists call domino-pieces {\it dimers}.
After writing [Z1], we realized that the same idea still applies 
to the problem of enumerating towers where instead of using dimers as pieces, one
uses trimes, 
or tetramers, or pentamers, etc.  Even more surprisingly, 
the [BeP] ideas may be used 
to enumerate towers using {\it any} given (finite) set of $k$-mers, where
{\it all} interfaces are allowed. By a $k$-mer we mean a $1 \times k$ rectangular piece,
where the case $k=2$ is a dimer (alias domino-piece).

This is all implemented in the Maple package {\tt TOWERS},
written by DZ, and linkable from the ``front'' of this article

{\eighttt http://www.math.rutgers.edu/\~{}zeilberg/mamarim/mamarimhtml/migdal.html} ,

that implements these beautiful human ideas. There you can also find lots of new and exciting
enumeration theorems generated by SBE, that readers are welcome to extend even further.

{\bf Background}

Recall that a domino (alias {\it dimer}) is a $1 \times 2$ rectangular piece.
Define an $i$-mer to be a $1 \times i$ rectangular piece.

Fix a positive integer $k$, and  suppose that we are allowed to use, as pieces, $i$-mers 
with $1 \leq i \leq k$. A {\bf tower} is a two-dimensional configuration where the bottom consists
of (one or more) {\bf contiguous} pieces (i.e., no gaps). 
A general tower is formed by starting with
a bottom floor (which is already called a tower), and, if desired,
building higher floors by placing non-overlapping pieces
on top of the currently highest floor, in such a way
that every piece in the newly created floor, 
must touch (at least) one piece of the floor right below it
(or else the piece would drop down).
Of course, the length of the common boundary between two touching
pieces from adjacent floors must be a strictly positive {\it integer}
$\leq k$.

On the computer, we describe a tower as a list of lists where on each floor, 
we list the locations $[x,x+i]$, of the $i$-mer whose left-end is at $x$.

For example, when $k=3$, the following is one such  (legal) tower:
$$
[[0,1],[1,2],[2,5],[5,7]],[[3,4], [6,9]],[[3,6]] \quad ,
$$
but
$$
[[1,2],[2,5],[5,7]],[[3,4], [6,9]],[[4,6]] \quad ,
$$
is {\bf not} legal since the piece $[4,6]$, on the third floor,  does not overlap with  any piece on the second floor.

Note that it is OK, with our current convention, for a piece to be perfectly aligned with a bottom piece. For example
$$
[[0,2],[0,2]] \quad ,
$$
is legal (but would not be in the ``xaviers'' counted by the original $3^n$ theorem).

Let the {\it weight} of one piece of size $i$ (i.e. of the form $[a,a+i]$, for some $a$), be
$ t^i z_i$, and let the  weight of a tower be the {\it product} of the weights of the individual pieces.

For example,
$$
weight \left ( [[1,2],[2,5],[5,7]],[[3,4], [6,9]],[[4,6]] \right )=
((t z_1 )(t^3 z_3)(t^2 z_2)) \cdot ((t z_1)(t^3 z_3)) \cdot (t^2 z_2))=
$$
$$
t^{12} z_1^2 z_2^2 z_3^2  \quad .
$$
We are interested in $M_k=M_k(t;z_1, \dots, z_k)$, the generating function (alias weight-enumerator) of the (``infinite'') set
of {\it all} legal towers formed from $i$-mers with $1 \leq i \leq k$. Then the coefficient of
$t^n$ would give us the {\it generating  polynomial}, for all towers with area $n$, and
setting all the $z_i$'s to $1$ we would get the number of towers of area $n$.

As in the B\'etr\'ema-Penaud approach (see [Z1]), let a {\it pyramid} be a tower whose
first floor only consists of {\it one} piece, and, let a {\it half-pyramid} be a pyramid 
none of whose floors have pieces that lie strictly to the left of the bottom piece.

Let $H=H(t; z_1, \dots, z_k)$, $P=P(t; z_1, \dots, z_k)$, and $M=M(t; z_1, \dots, z_k)$ be the weight-enumerators
of half-pyramids, pyramids, and towers, respectively. Then an almost {\bf verbatim}
(you do it!) argument shows that $H$ satisfies the following {\it algebraic equation}:
$$
H=\sum_{i=1}^{k} t^i z_i (1+H)^i \quad .
\eqno(HalfPyramids)
$$
Once we know $H$, we can get $P$ from:
$$
P=\frac{H}{1-\sum_{i=1}^k (i-1)t^i z_i (1+H)^i} \quad ,
\eqno(Pyramids)
$$
and finally
$$
M=\frac{P}{1-H} \quad .
\eqno(Towers)
$$

If the set of piece-sizes is not $\{1, \dots, k \}$, but an arbitrary finite set of positive integers, $S$, then of course
$$
H=\sum_{i \in S } t^i z_i (1+H)^i \quad .
$$
From $H$ we can get $P$:
$$
P=\frac{H}{1-\sum_{i \in S} (i-1)t^i z_i (1+H)^i} \quad,
$$
and finally
$$
M=\frac{P}{1-H} \quad .
$$
If we only want straight-enumeration, then set all  $z_i=1$, getting
$$
H=\sum_{i \in S } t^i(1+H)^i \quad .
\eqno(HalfPyramids1)
$$
From this we can get $P$:
$$
P=\frac{H}{1-\sum_{i \in S} (i-1)t^i (1+H)^i} \quad,
\eqno(Pyramids1)
$$
and finally
$$
M=\frac{P}{1-H} \quad .
\eqno(Towers1)
$$

It is very easy, for a computer, by iterating $H \rightarrow \sum_{i \in S } t^i(1+H)^i$, starting
with $H=0$, to crank out the first $200$ terms, but it gets slower and slower for
higher terms. But thanks to Comtet's famous theorem, implemented in the Salvy-Zimmermann Maple
package {\tt gfun} [SaZ], we can (rigorously) find a linear differential equation with polynomial
coefficients, and from this (still using {\tt gfun}) a linear recurrence equation 
for the enumerating sequence of half-pyramids,
that would give you, {\it much faster} than the algebraic equation, the $50000$-th term, say.

Also from $(HalfPyramids1)$, the algebraic equation for $H$, one can derive algebraic
equations for $P$ and $M$ using $(Pyramids1)$ and $(Towers1)$ respectively.

Furthermore, from these algebraic equations, one can use the beautiful methods described in
[FS] to derive asymptotics. Alternatively, one can derive them from the linear recurrences
by using the Maple package {\tt AsyRec} described in [Z2].

But an even more efficient approach is an {\bf empirical} one. 
First use the algebraic equations to crank out the first few terms of
the desired sequences, {\it guess}  linear recurrences, then use them to crank out {\bf many} terms, and use
empirical asymptotics to estimate the asymptotics. So everything is, if not ``rigorous'', at least
semi-rigorous, and easily {\it rigorizable}. Being {\it empiricists}, we prefer the latter method.

{\bf Encore: The One Piece-Size Case}

{\bf I. All Interfaces are allowed}

If there is only one piece, of size $k$, and {\it all interfaces are allowed} then the 
generating function, $H$, in the variable $b$, for the number of half-pyramids 
with $n$ pieces (so $b=t^k$) is:
$$
H=b (1+H)^k \quad ,
$$
and by {\it Lagrange Inversion}, we get the humanly-generated fact that
the number of half-pyramids (where all interfaces are allowed) using $n$ pieces,
each of length $k$, equals $\frac{(kn)!}{n!(kn-n+1)!}$, and
thanks to $P=\frac{H}{1-(k-1)H}$, we get  that the number of such pyramids is $\frac{1}{k} { {kn} \choose {n}}$.
For $k>2$, there is no `nice' expression for the number of towers with $n$ pieces,
but of course, 
using $M=\frac{H}{(1-(k-1)H))(1-H)}$, one can get  a linear recurrence,
either directly, or better still, empirically. 

For $k=2$ we have a {\bf nice surprise}, the number of domino towers with $n$ pieces,
where all interfaces are allowed, is $4^{n-1}$, in nice analogy with the fact that the
number of xaviers (where the exact-alignment interface is forbidden) is $3^{n-1}$.

{\bf II. All Interfaces are allowed, except Exact Alignment}

Recall that in the original, Viennot, scenario, that lead to the beautiful $3^{n-1}$ formula,
with dimers (dominoes), it was forbidden to place a piece exactly aligned with a piece on the
floor below it. If we impose this rule for the {\it one-piece case}, of size $k$,  (unfortunately, the
analysis is inapplicable for structures with more than one piece-size), we get the  equation
$$
H=b ((1+H)^k-H) \quad ,
$$
with, as above $P=\frac{H}{1-(k-1)H}$, $M=\frac{H}{(1-(k-1)H))(1-H)}$.

For $k>2$, things are no longer closed-form,  but 
we still get linear recurrences with polynomial coefficients, since the generating functions are algebraic, and hence $D$-finite.

{\bf References}

[BeP] J. B\'etr\'ema et J.-G. Penaud, 
{\it Mod\`eles avec particules dures, 
animaux dirig\'es et s\'eries en variables partiellement commutatives}
{\eighttt http://fr.arxiv.org/abs/math/0106210}

[Bo] Mireille Bousquet-M\'elou, 
{\it Rational and algebraic series in combinatorial enumeration}, 
invited paper for  the International Congress of Mathematicians 2006. 
Proceedings of the ICM. Session lectures, pp. 789--826. 
{\eighttt http://fr.arxiv.org/abs/0805.0588}.

[FS] P. Flajolet and R. Sedgewick, {\it ``Analytic Combinatorics''}, Cambridge University Press, 2009. \hfill \break
[freely(!) available on-line from {\eighttt http://algo.inria.fr/flajolet/Publications/book.pdf} ]

[GV] D. Gouyou-Beauchamps and G. Viennot, {\it Equivalence of the two dimensional directed animals
problem to a one-dimensional path problem}, Adv. in Appl. Math. {\bf 9} (1988), 334-357.

[SaZ] Bruno Salvy and Paul Zimmermann,
{\it Gfun: a Maple package for the manipulation of generating and holonomic functions in one variable},
ACM Transactions on Mathematical Software, {\bf 20} (1994), 163-177.

[V] G. X. Viennot, {\it Probl\`emes combinatoires pos\'es par la physique
statistique}, S\'eminaire N. Bourbaki, expos\'e $n^{o}$ 626, $36^{e}$ ann\'ee,
in Ast\'erisque $n^{o}$ 121-122 (1985) 225-246, SMF.
[available on-line]

[Z1] D. Zeilberger,
{\it The amazing $3^n$ theorem  and its even more Amazing Proof 
[discovered by Xavier Viennot and his \'Ecole Bordelaise gang]},
Personal Journal of Shalosh B. Ekhad and Doron Zeilberger.
\quad {\eighttt http://www.math.rutgers.edu/\~{}zeilberg/mamarim/mamarimhtml/bordelaise.html}

[Z2] Doron Zeilberger, 
{\it AsyRec: A Maple package for Computing the Asymptotics of Solutions of Linear Recurrence Equations with Polynomial Coefficients},
Personal Journal of Shalosh B. Ekhad and Doron Zeilberger.
\quad {\eighttt http://www.math.rutgers.edu/\~{}zeilberg/mamarim/mamarimhtml/asy.html}

\end